\documentclass{commat}

\usepackage[all]{xy}

\DeclareMathOperator{\Aut}{Aut}
\DeclareMathOperator{\Bir}{Bir}
\DeclareMathOperator{\Br}{Br}
\DeclareMathOperator{\Char}{char}
\newcommand{\cO}{\mathcal{O}}
\DeclareMathOperator{\cont}{cont}
\DeclareMathOperator{\dil}{dil}
\DeclareMathOperator{\Gal}{Gal}
\DeclareMathOperator{\Ker}{Ker}
\DeclareMathOperator{\Pic}{Pic}
\DeclareMathOperator{\Proj}{Proj}
\DeclareMathOperator{\res}{res}
\DeclareMathOperator{\Spec}{Spec}
\newcommand{\tto}{\longrightarrow}

\title{On birational automorphisms of Severi--Brauer surfaces}

\author{F.~V.~Weinstein}

\affiliation{N/A
    \email{felix.weinstein46@gmail.com}
    }

\abstract{
The generators of the group of birational automorphisms of any Severi\--Brauer surface non-isomorphic over an algebraically non-closed field to the projective plane are explicitly described.
    }



\DOI{10.46298/cm.9040}

\VOLUME{30}

\firstpage{1}

\begin{paper}

The classical theorem by M.~Noether states that the group of
birational automorphisms of the projective plane
\(\mathbf{P}_{\overline k}^2\) over an algebraically closed field \(\overline k\)
is generated by projective automorphisms and standard quadratic
Cremona transformations (see~\cite[Ch.~V, \S\S~5,6]{1}). In a~generalization of this
theorem to non-closed fields it is natural to consider, together
with projective plane, the \textit{Severi--Brauer surfaces}, i.e.,
surfaces defined over non-closed fields which become
isomorphic to the projective plane if their field of definition is
lifted up to its algebraic closure.

The arising situation has a~sapid cohomologic interpretation. Let
\(\overline G = \Gal(\overline k/k)\) be the Galois group of the
algebraic closure. Classes of Severi--Brauer surfaces over \(k\) can
be identified with elements of the
space \(H^1 (\overline G;PGL_3(\overline k))\), up to a~\(k\)-isomorphism. It is
known (see~\cite{3}) that the exact sequence of groups
\[
1\tto\overline k^\times\tto GL_3(\overline k)\tto PGL_3(\overline k)\tto 1
\]
induces an embedding
\(H^1 (\overline G;PGL_3(\overline k))\overset\delta\tto\Br(k)\),
where \(\Br(k) = H^2 (\overline G,\overline k^\times)\)
is the Brauer group of classes (up to an equivalence) of
central simple \(k\)-algebras. The image of the embedding \(\delta\)
consists exactly of the elements \(\gamma\in\Br(k)\) with Schur index (see~\cite{P}) which
divides \(3\), i.e., is equal to either 1 (the projective plane) or 3 (a Severi-Brauer surface).

Let \(S_2(\gamma)\) be the set of Severi--Brauer surfaces
corresponding to \(\gamma\), and let \(V\in S_2(\gamma)\).
The main result
of this paper is a~description of generators of the group of
birational automorphisms of \(V\) if \(\gamma\ne1\), i.e., of
Severi--Brauer surface non-isomorphic over \(k\) to the projective
plane. It is interesting that this description requires a~reference
to \(V'\in S_2(\gamma^{-1})\). The group of birational automorphisms
of \(V\) contains a~group of biregular automorphisms of \(V\) (described
by Theorem 3 of~\cite{3}) as a~subgroup.

This work had been carried out in 1970 as my MS Diploma thesis at the
Department of Mechanics and Mathematics, Moscow State University. D.~Leites translated it and
preprinted in proceedings of his ``Seminar on
Supersymmetries" (Reports of the Department of Mathematics, Stockholm University,
33/1989-2).

Recently the above-mentioned preprint of this text was cited
in an interesting paper by C.~Shramov~\cite{Shr}.
Since the result of this old work of mine is still useful,
I decided to update my preprint, make it available by putting it in arXiv,
and add to it comments I got meanwhile. First of all, from Torsten Ekedahl.

I wish to express my deep gratitude to my former
scientific advisor Prof.~Yu.\,I. Manin. I am very thankful to Torsten Ekedahl for his suggestions
how to simplify certain proofs in the above-mentioned preprint;
following his generous advice I cite his suggestions.
I am also thankful to D.~Leites for help
and to A.~Skorobogatov whose comment I got via D.~Leites.

\section{Birational maps of Severi--Brauer surfaces
corresponding to points of degree 3.}\label{s1}

Fix \(V\in S_2(\gamma)\) and \(V'\in S_2(\gamma^{-1})\). Theorem 2
of~\cite{3} and the theory of central simple algebras easily imply
that there always exist points of degree 3 on \(V\) and \(V'\), there
are no points of lesser degree, and degrees of all closed points
are multiples of 3.

The aim of this section is to associate with every point \(x\in V\) of
degree 3 a~birational map \(\varphi_3(x):V\tto V'\).
In what follows I always assume that the characteristics of the field \(k\) over
which \(V\) is defined is not 2 or 3. Thus, let \(x\) be a~point of
degree 3 belonging to \(V\). Let us perform a~monoidal transformation
\(\dil_x\colon V_1\tto V\) with the center at \(x\).

An exceptional curve of the first kind \(L\) which is different from \(\dil_x^{-1}(x)\)
belongs to \(V_1\). Contracting this curve leads again to a~Severi-Brauer surface,
and we thus obtain a~birational map
\(\beta_\alpha (x):V\tto V_\alpha '\).
By Theorem 5 of~\cite{3} either
\(V_\alpha '\in S_2(\gamma)\) or \(V_\alpha '\in S_2(\gamma^{-1})\).

\begin{lemma}\label{th1}
\(V_\alpha '\in S_2(\gamma^{-1})\).
\end{lemma}

\begin{proof}
Let \(k(x)\) be the field of quotients of the local ring at point \(x\) and \(K\)
the minimal normal extension of \(k\) containing \(k(x)\).
Let \(G = \Gal(K /k)\) and \(G_1 = \Gal(K/k(x))\). The natural homomorphism
\(H^2 (G;K^\times)\tto\Br(k)\) is an embedding whose image
contains \(\gamma\). Since it will not cause a~misunderstanding, the
preimage of \(\gamma\) will be also denoted by \(\gamma\).
Let \(\gamma\in H^2 (G;K^\times)\) be represented by a~system of factors
\(C_{\sigma,\tau}\in K^\times \), where \(\sigma, \tau\in G\). Let
\((K,G,\gamma)\) denote the \(k\)-algebra given by the system of relations
\[
\left\{
\begin{array}
{ll}
u_\sigma u_\tau = u_{\sigma\tau}\cdot C_{\sigma,\tau}&\\
cu_\sigma = u_\sigma\cdot c^\sigma, &\text{where \(c\in K^\times\)},
\end{array}
\right.
\]
i.e., the twisted product of \(K\) by \(G\) relative the system
\(C_{\sigma,\tau}\).

Let \(M\) be an irreducible right \((K,G,\gamma)\)-module.
By Lemma 2 in \S 5 of~\cite{R2} there exists \(T\in M\) such that
\(Tu_{\sigma_1} = T\cdot c_{\sigma_1}\) for all \(\sigma_1\in G_1\), where
\(c_{\sigma_1}\in K^\times \). Set
\(T^{u_\rho} = T\cdot u_\rho\); we see that
\(M = \oplus_{\rho\in G\mod G_1}T^{u_\rho}\cdot K\).

Since \(\Char k\ne3\), the field \(k (x)\) is separable over \(k\), and
hence the full system of representatives of classes \(G\mod G_1\) contains 3
elements: \(\rho_0\), \(\rho_1\) and \(\rho_2\). By setting
\(T_i := T^{u_{\rho_i}}\) we have \(M = T_0k\oplus T_1k\oplus T_2k\).

Let \(\mathbf P_K^2 = \Proj K[T_0,T_1,T_2]\), let \(\cO(2)\) be the
canonical invertible sheaf over \(\mathbf P_K^2 \) and let
\(\Gamma (\mathbf P_K^2,\cO(2))\) be the \(K\)-vector space of its global sections. The
subspace \(M'\subset\Gamma (\mathbf P_K^2,\cO(2))\), where
\[
M' = T_0T_1K \oplus T_1 T_2K\oplus T_0T_2K,
\]
defines a~linear system of conics on
the plane and the choice of the basis in it given by \(\{T_0T_1, \, T_1T_2, \, T_0T_2\}\) defines a~birational isomorphism
\({\overline\beta :\mathbf P_K^2 \tto\mathbf P_K^2}\)
which is the standard quadratic Cremona
transformation with the center in the triple
\[
\text{\(x_1 = (0:0:1)\),\
\(x_2 = (0:1:0)\),\ \(x_3 = (1:0:0)\).}
\]

Let \(L_i\) be the line in the plane \(\mathbf P_K^2 \)
containing 2 of these points different form \(x_i\) for \(i = 1,2,3\).
Evidently, the cycles \(C_ x = x_1 + x_2 + x_3\) and \(C_L = L_1+L_2+L_3\)
on the plane \(\mathbf P_K^2 \) are simple rational over \(k\) cycles with respect to the
Galois group of the covering \({\mathbf P_K^2 \cong V\otimes_kK\tto V}\).
Hence, the cycle \(C_x\) defines a~closed point on \(V\) which can
be considered coinciding with \(x\). Results of~\cite[\S~4]{R2} imply
that the birational map \(\overline\beta\) can be descended to a
birational map \(\beta :V\tto V_\alpha ''\), where
\(V_\alpha ''\in S_2(\gamma^{-1})\).

Due to the regularity of \(C_L\) it is c1ear that \(\varphi_3(x)\)
coincides, up to a~biregular isomorphism
\(V'_\alpha\tto V_\alpha ''\), with the map \(\beta_\alpha (x):V\tto V_\alpha '\)
constructed earlier.
\end{proof}

In March, 1985, A. Skorobogatov informed me of the following short proof of
Lemma~\ref{th1}; in particular, this enables one to get rid of
necessity to refer to~\cite{R2}.

By a~theorem of Manin (see~\cite[Ch.~IV]{Ma}), \(U = V\diagdown C_L\)
is a~principal homogeneous space over
2-dimensional torus \(T(K)\). Consider the following commutative
diagram:
\[
\xymatrix{
1\ar[r] &K^\times\ar[r] &GL_3(K)\ar[r] &PGL_3(K)\ar[r] &1\\
1\ar[r] &K^\times\ar[r]\ar[u] &R_{K/k(G_m(K))}\ar[r]\ar@{^{(}->}[u] &T(K)\ar[r]\ar@{^{(}->}[u] &1
}
\]
where \(R_{K/k}\) is the Weil functor, \(G_m(K)\) is the group of
\(K\)-points of the multiplicative group, see~\cite{5}. Passing to the
Galois cohomology we get, thanks to Hilbert's Satz 90,
the following diagram
\[
\xymatrix{
1\ar[r] &H^1 (G;PGL_3(K))\ar[rd]\\
  & &{\quad H^2 (G;K^\times)\subset \Br(k)}\\
1\ar[r] &H^1 (G;T(K)) \ar[uu]\ar[ru]
}
\]
The considered (Cremona) transformation acts on \(U\) and being lifted
to \(T(K)\) is the inversion \(x\mapsto x^{-1}\). Hence, in the group of
principal homogeneous spaces, \(H^1 (G;T(K))\), this transformation
induces an inversion, and therefore it does the same in the Brauer group.
Lemma~\ref{th1} is proved once more.

Now, take a~biregular isomorphism \(\psi_\alpha:V_\alpha '\tto V'\)
and define a~birational map
\[
\varphi_3^\alpha (x) = \psi_\alpha\circ\beta_\alpha (x):V\tto V'.
\]

The element \(\varphi_3^\alpha(x)\) of the set \(\Bir(V,V')\) of
birational maps is defined uniquely up to an action of the group
\(\Aut(V')\) of biregular automorphisms of \(V'\) on the set
\(\Bir(V,V')\).

Now, let \(Z(V) = Z(\gamma)\) be the group of cycles on \(V\)
(see~\cite{2}). It can be represented as a~direct sum
\(Z(\gamma) = \Pic(V)\oplus Z^0 (\gamma)\),
where \(Z_0(\gamma)\) is the group of
0-dimensional cycles on \(V\).

It will be convenient for us to
describe every 0-dimensional cycle on \(V\) as a~simple rational
cycle on
\(V\otimes_k\overline k\cong\mathbf P_{\overline k}^2 \).
The Picard group of \(V\) is isomorphic to a~free cyclic group
with the anticanonical class as a~generator, i.e.,
\(\Pic(V)\cong\mathbf Z(-\omega_\gamma)\).

\begin{lemma}\label{th2}
Let \(\varphi_3(x)_*:Z(\gamma)\tto Z(\gamma^{-1})\)
be the homomorphism induced by \(\varphi_3(x)\), let
\[
\alpha = -d\omega_\gamma -b(x_1+x_2+x_3)-
\sum\limits_{i\geqslant 4}b_ix_i\in Z(V).
\]
Then,
\[
\varphi_3(x)_*(\alpha) = -(2d-b)\omega_{\gamma^{-1}}-(3d-2b)(x_1'+x_2'+x_3')-
\sum\limits_{i\geqslant 4}b_ix_i',
\]
where \(x_1'+x_2'+x_3'\) is a~simple rational cycle on
\(V'\otimes_k\overline k\) and \({V'\in S_2(\gamma^{-1})}\); this curve is the
image of the exceptional curve of the first kind under \(\cont_L\).
\end{lemma}

\begin{proof}
It is subject to a~simple calculation, see~\cite{1}.
\end{proof}

\section{Birational maps of Severi--Brauer surfaces corresponding
to points of degree 6}

To every point \(x\in V\) of degree 6 a birational map \(\varphi_6(x):V\tto V'\) can also be assigned, where \(V'\) is a~fixed element of \(S_2(\gamma^{-1})\). We will need the following statement.

\begin{lemma}\label{th3}
Consider the fiber product
\(V\otimes_k\overline k\cong\mathbf P_{\overline k}^2 \)
and its projection onto the first factor
\(p:\mathbf P_{\overline k}^2 \to V\). Let \(x\in V\) be a~closed point of degree \(6\) and define \(p^{-1}(x) := (x_0,\ldots,x_5)\), where the \(x_i\)
are closed points in \(\mathbf P_{\overline k}^2 \). Then, no \(3\) points \(x_i\) belong to one line and all \(6\) points do not belong to a~conic.
\end{lemma}

\begin{proof}
If no 3 points belong to one line and all 6 belong to
a conic, then these points uniquely define this conic. This conic
defines, on \(V\), a~simple rational cycle over \(k\). The divisor
corresponding to this cycle is of degree 2 contradicting
Proposition 13 of~\cite{3}.

Now suppose that there is a~line
in \(\mathbf P_{\overline k}^2 \) containing 3 points \(x_i\).
The totality of all lines with
pair-wise distinct points \(x_i\) forms a~cycle rational over \(k\). The
degree of the divisor corresponding to this cycle is equal to the
number of these lines and due to~\cite[Proposition 13]{3} should
be a~multiple of 3. Since the number of lines which connect
pair-wise distinct points \(x_i\) does not depend on the order of
these points on lines, it follows that from the combinatorial
point of view only the following \(9\) cases are possible:
\begin{enumerate}
\item[1)] All \(6\) points belong to one line.

\item[2)] There is a~line containing exactly \(5\) points.

\item[3)] There is a~line containing exactly \(4\) points, and there is no line containing \(3\) points.

\item[4)] There exists a~line containing exactly \(4\) points, and a~line containing exactly \(3\) points.

\item[5)] There exist exactly two lines each of them containing \(3\) points,
together they contain \(5\) points, and there are no lines containing
more than \(3\) points different from these two lines.

\item[5+i)] There exist exactly \(i\) lines, where \(i = 1,2,3\) or \(4\), each of
them containing exactly \(3\) points, there are no lines containing more
points, and the case 5) fails.
\end{enumerate}

Let \(n_j\) be the number of lines connecting pairs of points of the set
in case \(j)\). Then, it is easy to see that the following relations
hold:
\begin{center}
\begin{tabular}{crrrrrrrrr}
\(j\)&1&2&3&4&5&6&7&8&9\\
\hline \(n_j\)&1&6&10&8&11&13&11&9&7
\end{tabular}
\end{center}

Due to the above, it follows that only cases 2) and 8) can hold.

Case 2). Let \(x_1\), \(\ldots \), \(x_5\) be points belonging to the line
\(L\). By transitivity of the \(\Gal (\overline k/k)\)-action on
\((x_1,\ldots,x_6)\) we can find \(g\in\Gal(\overline k/k)\) such
that \(g(x_1) = x_6\). Then, \(g\) transforms \(L\) into a~line passing
through \(x_6\). But \(L\) contains 5 points and not all 6 points belong
to one line, hence a~contradiction.

Case 8). Let \(L_1\), \(L_2\), and \(L_3\) be lines passing through \(x_1\),
\(x_2\), \(x_3\), through \(x_3\), \(x_4\), \(x_5\), and through \(x_5\), \(x_6\), \(x_1\),
respectively.
We have \(g(L_i) = L_j\), where \(i\), \(j = 1,2,3\),
for any \(g\in\Gal (\overline k/k)\). Moreover,
\(g(L_i\cap L_j) = L_{i^\prime}\cap L_{j^\prime}\) since if it is not so, then
through the point not representable in the form \(L_i\cap L_j\) two
lines from the set \(\{L_1,L_2,L_3\}\) pass, but then
\(\Gal(\overline k/k)\) does not act transitively on \(\{x_1,\ldots,x_6\}\)
contradicting the simplicity of the cycle \(x_1+\ldots+x_6\).
\end{proof}

By Lemma~\ref{th3} to any point \(x\in V\) such that \(\deg_k(x) = 6\) we
can assign a~birational map \(\varphi_6(x):V\tto V'\).

Let \(\dil_x:\overline V\tto\mathbf P_{\overline k}^2 \) be a~monoidal
transformation with the center at the points \(x_1, x_2, \dotsc, x_6\),
and let \(Q_i\) be a~conic in
\(\mathbf P_{\overline k}^2 \) containing
\(\{x_1,\ldots,x_{i-1},x_{i+1},\ldots,x_6\}\).
Let \(S_i\) be the proper preimage of \(Q_i\) with respect to \(\dil_x\).
The curve \(S = \underset {0\leqslant i\leqslant 6}\cup S_i\) is contractible
since it is an exceptional curve of the first kind.
Let \(\cont_S:V\tto\mathbf P_{k}^2 \) be a~contraction morphism of \(S\) and
\(\overline x_i = \cont_S(S_i)\). Thus, we have a~birational isomorphism
\[
\overline\beta = \cont_S\cdot\dil_x^{-1}:\mathbf P_{\overline
k}^2 \tto\mathbf P_{\overline k}^2.
\]
It is easy to compute the value of the homomorphism
\({\overline\varphi_6(x)_*:Z(\mathbf P_{\overline k}^2)\tto Z(\mathbf P_{\overline k}^2)}\)
at the anti\-canonical class
\(-\omega_{\mathbf P_{\overline k}^2} = -\omega\):
\begin{equation}
\overline\beta_*(-\omega) = -5\omega -6(\overline x_1+\ldots
+\overline x_6).\label{eq1}
\end{equation}
Since \(Q_1 + \dotsc + Q_6\) is a~simple and rational cycle over \(k\),
then so is \(S_1 + \dotsb + S_6\). Hence, \(\overline\varphi_6(x)\) can
descend to a~birational isomorphism \(\beta_\alpha (x):V\tto
V_\alpha '\) which is a~composition of the blowing up of the point
\(x\in V\) such that \(\deg_k(x) = 6\) and a~contraction of the
exceptional curve \(S\) of the first kind.
By Theorem 5 of~\cite{3}
either \(V_\alpha '\in S_2(\gamma)\) or \(V_\alpha '\in S_2(\gamma^{-1})\).

\begin{lemma}\label{th4}
\(V_\alpha '\in S_2(\gamma^{-1})\).
\end{lemma}

\begin{proof}
Let \(k(x)\) be the field of quotients of the local ring at point \(x\). Consider
the minimal normal extension \(K\) of the field \(k\) such that
\(k\subset k(x)\subset K\). Since \([k(x):k] = 6\), it follows that
\(\Gal(K/k)\) contains at least one Sylow 3-subgroup.
Let \(G(3)\subset\Gal(\overline k/k)\) be a~Sylow 3-subgroup, and \(K_3\subset K\) the
subfield of elements fixed under \(G(3)\).

Set \(G = \Gal(\overline k/k)\) and \(G_1 = \Gal(\overline k/K_3)\).
Consider the restriction homomorphism
\[
\res:H^2 (G;\overline k^\times)\tto H^2 (G_1;\overline k^\times).
\]
Let
\[
\res(\gamma) = \gamma_1.
\]
Since \([K_3:k]\not\equiv 0\mod3\), then \(\gamma_1\ne1\). Setting
\(V_1 = V\otimes_kK_3\) we see that \(V_1\in S_2(\gamma_1)\). The fiber
over \(x\) of the projection \(p_1:V_1\tto V\) is isomorphic to
\[
\Spec(K_3)\underset {\Spec(k)}\times \Spec(k(x)) = \Spec(K_3\otimes_kk(x)).
\]
It follows from the construction that
\[
\Spec(K_3\otimes_kk(x)) = \Spec(K_1)\amalg\Spec(K_2),
\]
where \([K_i:K_3] = 3\) for \(i = 1,2\). Thus, \(p^{-1}(x_0) = \{y_1,y_2\}\)
and \(\deg_{K_3}(y_i) = 3\) for \(i = 1,2\).
Since the diagram
\[
\xymatrix{
&\mathbf{P}^2_{\overline{k}}\ar[ld]_p\ar[rd]^p \\
V_2\ar[rr]^{p_1} & & V_2
}
\]
is commutative, we can assume that
\[
\text{\(p^{-1}(y_1) = \{x_1,x_2,x_3\}\) and
\({p^{-1}(y_2) = \{x_4,x_5,x_6\}}\).}
\]

Now suppose Lemma~\ref{th4} fails. Then, \(V_\alpha '\in S_2(\gamma)\) and
we have a~birational isomorphism
\[
\beta_\alpha^1 = \beta_a(x)\otimes_kK_3:V_1\tto X_1,
\]
where \(X_1\in S_2(\gamma_1)\). By eq.\eqref{eq1}
\[
\beta_{\alpha *}^1 (-\omega_{\gamma_1}) = -5\omega_{\gamma_1}-
6(\overline x_1+\ldots +\overline x_6).
\]
For \(X'_1\in S_2(\gamma_1^{-1})\) and \(\overline y_0 = \{\overline x_1,\overline x_2,\overline x_3\}\in X_1\), if follows from Lemma~\ref{th2} and the birational isomorphism \({\psi_3(\overline y_0):X_1\to X_1'}\) that
\[
\varphi_3(\overline y_0)_*\circ\beta_{\alpha *}'(-\omega
_{\gamma_1}) = -4\omega_{\gamma_1^{-1}}-3(\overline
x_1'+\overline x_2'+\overline x_3')- 6(\overline x_4'+\overline
x_5'+\overline x_6').
\]
Then, for the birational isomorphism
\({\varphi_3(\overline y_0'):X_1'\to X_1}\), where \(\overline y_0'\! = \! \{\overline x_4', \overline x_5', \overline x_6'\}\), we have
\[
\varphi_3(\overline y_0')_*\circ
\varphi_3(\overline y_0)_*\circ\beta_{\alpha *}'(-\omega
_{\gamma_1}) = -2\omega_{\gamma_1}-3(\overline x_1''+\overline
x_2''+\overline x_3'').
\]
Finally, for the birational isomorphism
\(\varphi_3(\overline y_0''):X_1\to X_1'\), where
\(y_0'' = \{\overline x_1'',\overline x_2'',\overline x_3''\}\), we have
\[
\varphi_3(\overline y_0'')\circ\varphi_3(\overline y_0')_*\circ
\varphi_3(\overline y_0)_*\circ\beta_{\alpha *}^1 (-\omega
_{\gamma_1}) = -\omega_{\gamma_1^{-1}}.
\]
Thus, we got a~biregular isomorphism \(V_1\tto X_1'\), where
\(X_1'\in S_2(\gamma_1^{-1})\). This is a~contradiction.
\end{proof}

As it had been done with points of degree 3, we will associate with
every \(x\in V\) a~birational map \(\varphi_6(x):V\tto V'\).

\begin{lemma}\label{th5} Let
\(\varphi_6(x)_*:Z(\gamma)\tto Z(\gamma^{-1})\)
be the homomorphism induced by \(\varphi_6(x)\).
Then,
\begin{multline*}
\varphi_6(x)_*\Big(-d\omega_\gamma -b(x_1+\ldots +x_6)-
\sum\limits_{i\geqslant 7}b_ix_i\Big)\\
= -(5d-4b)\omega_{\gamma
^{-1}}-(6d-5b)(x_1'+\ldots +x_6')-\sum\limits_{i\geqslant
7}b_ix_i',
\end{multline*}
where \(s = x_1'+\ldots+x_6'\) is a~simple rational over \(k\) cycle on
\(V'\otimes_k k'\) with \(V'\in S_2(\gamma^{-1})\) which is the
image of contraction of the exceptional curve of the first kind defined by \(s\).
\end{lemma}

\begin{proof}
It suffices to represent the map
\(\overline\varphi_6(x):\mathbf P_{\overline k}^2 \tto\mathbf P_{\overline k}^2 \)
as a~product of 3 quadratic Cremona transformations.
\end{proof}

\section{Proof of the main theorem}
Let \(V\in S_2(\gamma)\) and \(V'\in S_2(\gamma^{-1})\) be fixed
surfaces. The said earlier can be summed up as follows: To every
point \(x\in V\) such that \({\deg x = 3}\) or 6 there uniquely
corresponds an orbit of the left action of the group
\(\Aut(V')\) of biregular automorphisms on the set \(\Bir(V,V')\) of birational maps
\(V\tto V'\). For every such point, we will choose, once and for all,
an element of the corresponding orbit. It is a~birational map
\(\varphi_i(x):V\tto V'\), where \(i = \deg_kx\). Similarly, for any
\(y\in V'\), we will construct a~map \(\varphi_j(y):V'\tto V\).

Our goal is the proof of the following statement.

\begin{theorem}\label{th6} The group \(\Bir (V)\) of birational
automorphisms of \(V\) is generated by its subgroup of biregular
automorphisms \(\Aut(V)\), and automorphisms \(\varphi_i(x)\) and
\(\varphi_i(y)\) for all \(x\in V\) and \(y\in V'\) described above.
\end{theorem}

The proof of this theorem follows from a~series of lemmas. Let
\({g\in\Bir(V)}\) be an arbitrary birational automorphism and let the value of the
homomorphism \(f_*:Z(\gamma)\to Z(\gamma)\) at \(-\omega_\gamma\) be equal to \(-d\omega_\gamma-\sum\limits b_ix_i\).

\begin{lemma}\label{th7}
The following relations hold:
\begin{align}
 &9d^2 -\sum\limits b_i^2 = 9;\label{eq2}\\
 &9d-\sum\limits b_i = 9.\label{eq3}
\end{align}
\end{lemma}

\begin{proof}
Let us use the fact that \(f_*\) preserves
the arithmetic genus and the index of intersection of cycles on
\(V\), see~\cite{2}. We have
\[
{\left(-d\omega_\gamma -\sum\limits b_ix_i\right)}^2
= 9d^2 -\sum\limits b_i^2 = {(-\omega_\gamma)}^2 = 9
\]
yielding relation~\eqref{eq2}. Further,
\[
p_a(-d\omega_\gamma -\sum\limits
b_ix_i) = \frac92(d^2 -d)+1-\frac12\sum\limits b_i(b_i-1) = p_a(-\omega
_\gamma) = 1
\]
or \(9d^2 -\sum\limits b_i^2 = 9d-\sum\limits b_i\). Taking relation~\eqref{eq2}
into account we get relation~\eqref{eq3}.
\end{proof}

\begin{lemma}\label{th8} Let \(b = \max b_i\). Then,
\begin{equation}
b\geqslant d+1
\label{eq4}
\end{equation}
\end{lemma}

\begin{proof}
Indeed, relation~\eqref{eq2} implies \(9d^2 -b\sum\limits b_i\leqslant 9\).
Taking relation~\eqref{eq3} into account we see that \(9d^2 -9\leqslant b(9d-9)\) yielding inequality~\eqref{eq4}.
\end{proof}

\begin{lemma}\label{th9}
Let \(x_{i_0}\) be a~point of the cycle
\(f_*(-\omega_\gamma)\) whose coefficient is \(b\), i.e., \(x_{i_0}\) is
a point of maximal multiplicity. Then, \({\deg_kx_{i_0}<9}\).
\end{lemma}

\begin{proof}
Let \(n = \deg_kx_{i_0}\). Then, relation~\eqref{eq3} can be rewritten as \(9d-nb-\sum\limits_{i\ne i_0}b_i = 9\). Thanks to inequality~\eqref{eq4} we get
\[
9d-n(d+1)-\sum\limits_{i\ne i_0}b_i\geq9\text{\quad or \quad }-\sum\limits_{i\ne
i_0}b_i\geqslant 9-9d-n(d+1).
\]
Setting \(n = 9+3l\) we rewrite the
latter expression in the form
\[
-\sum\limits_{i\ne i_0}b_i\geqslant 18+3l(d+1)
\]
which is false if \(l\geqslant 0\) since
\(b_i\geqslant 0\), see~\cite [Corollary 1.18]{2}. Hence, \(l<0\), and then
\(\deg_kx_{i_0} = n\leqslant 6\), as was required.
\end{proof}

\noindent\noindent%
\textit{Proof of Theorem~\ref{th6}. }
It is well known (\cite{1}) that the point \(x\) of
maximal multiplicity of the cycle \(g_*(-\omega_\gamma)\)
belongs to \(V\). By Lemma~\ref{th9} the degree of this point is equal to
either 3 or 6. Therefore, applying the homomorphism \(\varphi_i(x)\)
to \(g_*(-\omega_\gamma)\) we will diminish, thanks to
Lemmas~\ref{th2},~\ref{th5},~\ref{th8}, the absolute value of the
coefficient of \(\omega_{\gamma^{-1}}\). Repeatedly applying this
procedure we will diminish the degree \(d\) to~1. Then, \(b_i = 0\) for all \(i\), as follows
from relation~\eqref{eq2}. Finally, we get either
\begin{equation}
\left\{\prod\limits f_{j_k,i_k}(y,x)_*\right\}\circ
g_*(-\omega_\gamma) = -\omega_\gamma
\label{eq5}
\end{equation}
or
\begin{equation}
\left\{\varphi_i(x)\circ\prod\limits f_{j_k,i_k}(y,x)_*\right\}
\circ g_*(-\omega_\gamma) = -\omega_\gamma.\label{eq6}
\end{equation}

Formula~\eqref{eq6} leads to a~contradiction since \(V\) and
\(V'\in S_2(\gamma^{-1})\) are not biregularly isomorphic.

Formula~\eqref{eq5}
implies that
\(\left\{\prod\limits f_{j_k,i_k}(y,x)\right\}\circ f\in\Aut(V)\).
Applying transformations inverse to \(f_{j_k,i_k}(y,x)\)
and isomorphisms \(\varphi_i\) to the left-hand side, we get the
required.
\hfill\qed

\begin{remark}
On automorphisms of similar (Del Pezzo) surfaces, see~\cite{6}.
\end{remark}

\section{Appendix. T.~Ekedahl's comments}

Most part of your paper is devoted to proof of Lemmas~\ref{th1} and~\ref{th4}.
The proof is overcomplicated with a~long and rather ugly
division into cases in the proof of Lemma~\ref{th3}.
It is possible to give a~short, uniform and conceptual treatment of these lemmas.

The main idea is to exploit the simple fact the Picard group
is of rank 2 for any \(k\)-surface \(X\) obtained by blowing up
a closed point on a~Severi--Brauer surface \(V\).
If we first look at Lemma~\ref{th3}, then \(X\)
is obtained blowing up a~closed point of order \(6\).

The statement of the lemma is equivalent to, and may be
replaced by, any of the statements:
\begin{enumerate}
\item[a)] \(-K_X\) is ample,
\item[b)] \(-K_X\) is very ample,
\item[c)] \(X\) is isomorphic to a cubic surface (cf. the English version of Manin's book~\cite[Ch.~IV, \S~24]{Ma}).
\end{enumerate}

It is easy to verify the statement a)~by using the ampleness criterion of
Moishezon and Nakai
(cf. proof of Statement 24.5.2 in op. cit. or p.~365
in Hartshorne's book~\cite{H}) and the fact that \(\rm{rk}\Pic X = 2\).

It is well known (cf. the comment on ``Sch\"afli's double-six''
in Hartshorne's book~\cite{H}) that there is a~natural set of
six (conjugated) lines on \(X\) complementing the six exceptional
\(\overline k\)-lines of \(X \longrightarrow V\).
In Lemma~\ref{th4} you study the surface \(V'\) obtained by contracting
these complementary lines. This is (cf. Manin's book~\cite{Ma})
a del Pezzo \(k\)-surface of degree~9, i.e., a~Severi--Brauer \(k\)-surface.
We have, therefore, two elements \(\{V\}\) and
\(\{V'\}\) in \(H_\text{\'et}(k,PGL_3)\) corresponding to the \(k\)-isomorphism
classes of \(V\) and \(V'\) (cf. Milne's book~\cite[p.~134]{Mi}).

In Lemma~\ref{th4} you consider the images
\(\langle V \rangle \) and \(\langle V' \rangle \) of the elements \(\{V\}\)
and \(\{V'\}\) in \(H^2_\text{\'et}(K,G_m)\) in the cohomology sequence
corresponding to the exact sequence of \'etale sheaves
(cf.~\cite[p.~142]{Mi}):
\[
1\tto G_m\tto GL_3\tto PGL_3\tto 1.
\]
You prove that these images are inverses of each other.

It is possible to give a~much more natural approach and prove
Lemmas~\ref{th1} and~\ref{th4} at the same time.
To begin with, the definition of the Brauer--Severi
schemes implies that \(\langle V \rangle \) (resp. \(\langle V' \rangle \))
belongs to the kernel of the natural map from
\(H_\text{\'et}^2 (K,G_m)\) to \(H_\text{\'et}^2 (V,G_m)\)
(resp. to \(H_\text{\'et}^2 (V',G_m)\)): use the fact
that the image of \(\{V\}\) in \(H_\text{\'et}^1 (V,PGL_3)\)
comes from \(H_\text{\'et}^1 (V,GL_3)\).
But it is known (\cite[p.~106]{Mi}) that \(H_\text{\'et}^2 (V,G_m)\) injects into \(H_\text{\'et}^2 (K(V),G_m) = H_\text{\'et}^2 (\overline k(V'),G_m)\)
and this implies that
\(\langle V' \rangle\in\Ker(H_\text{\'et}^2 (k,G_m)\tto H_\text{\'et}^2 (V,G_m))\).

The Hochschild--Serre spectral sequence
\[
H_\text{\'et}^p (k,H_\text{\'et}^q (\overline V,G_m))\Rightarrow H_\text{\'et}^{p+q}(V,G_m),
\]
see~\cite[p.~105]{Mi}, yields that
\[
\Ker(H_\text{\'et}^2 (k,G_m) \tto H_\text{\'et}^2 (V,G_m)) = \mathbb Z/3\mathbb Z,
\]
and since \(H_\text{\'et}^1 (k,GL_3) = 1\)
(cf.~\cite[p.~124]{Mi}) and \(V(k) = V'(k) = \emptyset\), it follows
that \(\langle V\rangle\) and \(\langle V'\rangle\) are non-trivial in \(H_\text{\'et}^2 (k,G_m)\).

It suffices to prove that \(\langle V\rangle\ne\langle V'\rangle\)
in \(H_\text{\'et}^2 (k,G_m)\). To show this, use the following
commutative diagram of \'etale sheaves over \(V\), \(V'\) and \(X\) (cf.~\cite[p.~143]{Mi}):
\begin{equation}\label{eq_last}
\xymatrix@R8mm{
&1\ar[d]&1\ar[d]\\
1\ar[r]&{\mu_3}\ar[d]\ar[r]&SL_3\ar[d]\ar[r]&PGL_3\ar[d]^{id}\ar[r]&1\\
1\ar[r]&G_m\ar[d]^3 \ar[r]&GL_3\ar[d]^{\det}\ar[r]&PGL_3\ar[r]&1\\
&G_m\ar[d]\ar@{ = }[r]&G_m\ar[d]\\
&1&1
}
\end{equation}

We then consider the image of \(\{V\}\in H_\text{\'et}^1 (k,PGL_3)\) under the composite map
\[
H_\text{\'et}^1 (k,PGL_3) \tto H_\text{\'et}^1 (V,PGL_3) \tto H_\text{\'et}^2 (V,\mu _3).
\]
This image lies in the kernel of the map
\(H_\text{\'et}^2 (V,\mu _3)\tto H_\text{\'et}^2 (V,G_m)\)
since the pullback of \(\{V\}\) in \(H_\text{\'et}^1 (V,PGL_3)\)
comes from an element in \(H_\text{\'et}^1 (V,GL_3)\).
We have, therefore, by the first column a~well-defined element
\([D_\gamma ]\) in \(\Pic V/3\Pic V\) corresponding to the
image of \(\{V\}\) in \(\ker (H_\text{\'et}^2 (V,\mu _3)\tto H_\text{\'et}^2 (V,G_m))\)
and we obtain in the same way an element
\([D_{\gamma '}]\) in \(\Pic V'/3\Pic V'\) from \(\{V'\}\).
If \(\{V\} = \{V'\}\) in \(H_\text{\'et}^1 (k,PGL_3)\),
then the pullbacks of \([D_\gamma ]\) and \([D_{\gamma '}]\) in \(\Pic X/3\Pic X\) must coincide.

But it is easy to compute \([D_\gamma ]\) and \([D_{\gamma '}]\).
We already noted that the image of
\(\{V\}\) in \(H_\text{\'et}^1 (V,PGL_3)\) comes from
\(H_\text{\'et}^1 (V,GL_3) = H_\text{zar}^1 (V,GL_3)\)
and it is known (cf., e.g., the end of Quillen's article~\cite{Q})
that one may choose an element in \(H_\text{zar}^1 (V,GL_3)\)
corresponding to a~vector bundle \(J_\gamma \) coming from a~natural extension
\[
0\tto\Omega_\gamma\tto J_\gamma\tto\cO_\gamma\tto 0,
\]
where \(\Omega_\gamma\) is the cotangent bundle of \(V\).
This implies by the diagram~\eqref{eq_last}
that \([D_\gamma ]\) is equal to the class of the line bundle
\(\det(J_\gamma)\), i.e., to the image
\([-K_\gamma]\in \Pic V/3\Pic V\) of the anti-canonical
class \(-K_\gamma \). But it is easy to see that the images of
\([-K_\gamma]\) and \([-K_{\gamma '}]\) in \(\Pic X/3\Pic X\)
do not coincide (here it is essential to consider \(\Pic X/3 \Pic X\) and not
\(\Pic\overline X/3\Pic\overline X\)).
This proves that \(\{V\}\ne\{V'\}\), thereby completing the proof of
Lemma~\ref{th1}. The same proof works for
Lemma~\ref{th4}; the only difference being that \(X\) is of degree~6.

It would also be useful for the reader if you include more modern references.
It would be valuable to have a~reference to the recent survey
article by Manin and Tsfasman~\cite{MTs}, so that the reader
can compare with other papers about birational automorphisms
on rational varieties like the ones by Iskovskikh and Manin (see also more recent papers~\cite{DI},~\cite{ISh}).

It would also be useful if you include a~reference to the excellent
article by M.~Artin~\cite{Ar} (and also to the one by Amitsur~\cite{Am})
about Severi--Brauer varieties in Springer LNM vol.~917.
Finally, I recommend you to make a~fuller use of
\'etale cohomology (cf. Milne's book~\cite{Mi})
which is the natural language for many of the results and arguments of the paper.

\EditInfo{%
    August 06, 2019}{%
    March 18, 2021}{%
    Dimitry Leites
}

\end{paper}